\documentclass{article}

\usepackage [preprint] {neurips_2026}

\usepackage[utf8]{inputenc} 
\usepackage[T1]{fontenc}    
\usepackage{hyperref}       
\usepackage{url}            
\usepackage{booktabs}       
\usepackage{amsfonts}       
\usepackage{nicefrac}       
\usepackage{microtype}      
\usepackage{xcolor}         
\usepackage{amsmath}
\usepackage{amsthm}
\usepackage{graphicx}
\usepackage[ruled,vlined,linesnumbered]{algorithm2e}
\usepackage{subcaption}
\usepackage{booktabs}
\usepackage{multirow}
\usepackage{makecell}

\theoremstyle{plain}

\theoremstyle{definition}

\theoremstyle{remark}

\title{Learning to Cut: Reinforcement Learning for Benders Decomposition}

\author{%
  Haochen Cai \\
  Department of Integrated Systems Engineering \\
  The Ohio State University \\
  Columbus, OH, USA, 43210 \\
  \texttt{cai.1256@osu.edu}
  \And
  Xian Yu \\
  Department of Integrated Systems Engineering \\
  The Ohio State University \\
  Columbus, OH, USA, 43210 \\
  \texttt{yu.3610@osu.edu}
}

\begin{document}

\maketitle

\begin{abstract}
Benders decomposition (BD) is a widely used solution approach for solving two-stage stochastic programs arising in real-world decision-making under uncertainty. However, it often suffers from slow convergence as the master problem grows with an increasing number of cuts. 
In this paper, we propose Reinforcement Learning for BD (RLBD), a framework that adaptively selects cuts using a neural network-based stochastic policy. The policy is trained using a policy gradient method via the REINFORCE algorithm. We evaluate the proposed approach on a two-stage stochastic electric vehicle charging station location problem and compare it with vanilla BD and LearnBD \citep{jia2021benders}, a supervised learning approach that classifies cuts using a support vector machine. Numerical results demonstrate that RLBD achieves substantial improvements in computational efficiency and exhibits strong generalization to problems with similar structures but varying data inputs and decision variable dimensions.
\end{abstract}

\section{Introduction}
Two-stage stochastic programs provide a general framework for modeling decision-making under uncertainty, where decisions are separated into \textit{here-and-now} (first-stage) decisions, made before uncertainty is realized, and \textit{wait-and-see} (second-stage) decisions, which adapt to observed outcomes through recourse actions after the uncertainty unfolds. This modeling paradigm has been widely used in a broad range of applications, including facility location problems \citep{magnanti1981accelerated,fischetti2016benders,fischetti2017redesigning}, supply chain management \citep{keyvanshokooh2016hybrid,basciftci2023resource}, scheduling and routing problems \citep{adulyasak2015benders,bodur2017mixed,yu2021integrated,yu2023time}, energy system operations \citep{alguacil2000multiperiod,wang2011chance,daneshvar2020two}, among others. A key challenge in these problems is that the second-stage decisions depend on random realizations, leading to large-scale optimization models. As the number of scenarios increases, directly solving the resulting extensive-form problem becomes computationally prohibitive. 

Benders decomposition (BD) has emerged as a popular method for solving two-stage stochastic programs by decomposing the original problem into a master problem and a set of subproblems \citep{benders1962partitioning}. Specifically, at each iteration, the algorithm generates cuts that approximate the recourse function (the second-stage value function) and iteratively tighten the master problem. As a result, the computational performance of BD depends critically on the quality of the cuts that are added during this iterative process. In the multi-cut formulation, the algorithm generates separate cuts for each scenario or subproblem, producing a tighter and more accurate representation of the recourse function; however, adding all of them can significantly increase the size of the master problem and slow down subsequent solves. On the other hand, in the single-cut formulation, all subproblem (scenario-wise) information is aggregated into one cut per iteration, leading to a smaller master problem but potentially weaker approximations of the recourse function. This introduces a trade-off between strengthening the master problem and maintaining computational efficiency. We refer to \cite{birge2011introduction} for a detailed introduction of these two formulations.

In this work, we model cut selection in BD as a sequential decision-making problem. At each iteration, the state encodes current convergence statistics, along with features of the master problem, subproblems, and individual cuts, while the action corresponds to selecting a subset of candidate cuts to add to the master problem. Based on this formulation, we develop a reinforcement learning (RL) approach that trains an optimal cut-selection policy using policy gradient algorithms. Our contributions can be summarized as follows:
\begin{itemize}

    \item \textbf{Stochastic cut-selection policy. } Different than the existing literature in BD that selects cuts deterministically, we consider a stochastic policy that selects each cut with a certain probability,  based on global algorithmic information and local cut features. 
    This allows the policy to exploit the most informative cuts given the current information while still maintaining the ability to explore alternative cuts that may become important later.

    \item \textbf{Strong generalization to problems with similar structures.} Our numerical results demonstrate that the trained cut-selection policy exhibits strong generalization to problems with similar structures but different data inputs and decision variable dimensions. 
    
    \item \textbf{Substantial improvements over existing benchmarks.} Compared to existing benchmarks such as vanilla BD and LearnBD \citep{jia2021benders}, our RLBD reduces computational time by up to a factor of five on most testing instances and achieves smaller optimality gaps when all methods fail to converge within the prescribed time limit.
\end{itemize}


\section{Related Work}
A variety of classical techniques have been developed to accelerate the convergence of BD. Among them, cut management is a key technique that carefully controls which cuts to add and retain in the master problem. Instead of adding all candidate cuts from solving the subproblems, many approaches add only a selected set of informative cuts at each iteration, such as the most violated cuts or the Pareto-optimal cuts, to improve convergence without excessively enlarging the master problem \citep{magnanti1981accelerated, fischetti2017note}. Another widely used heuristic is to remove redundant or inactive cuts that no longer contribute to tightening the formulation, thereby controlling the size of the master problem. On the other hand, scenario grouping techniques have been studied in the literature to improve computational efficiency by aggregating multiple cuts from similar scenarios \citep{adulyasak2015benders,maheo2019benders,ramirez2023benders}. An alternative stream of research focuses on reducing the solution time per iteration by solving the master problem or subproblem inexactly, leading to the so-called inexact BD \citep{geoffrion1974multicommodity,easwaran2009tabu,poojari2009improving}.


Recent works have explored integrating machine learning (ML) and RL into BD to improve its efficiency. \citet{jia2021benders} proposed a learning-enhanced framework (LearnBD) that trains a support vector machine to classify cuts and selectively adds those likely to be valuable to the master problem, thereby reducing the size of the master problem. \citet{choi2025learning} replaced the exact subproblem solver with a pretrained transformer model, accelerating the convergence of BD by generating high-quality approximate solutions for the subproblems. 
Different than these ML-based prediction models, \cite{li2025rl} applied RL to learn policies to adaptively select optimality gap tolerance for the master problem based on the features of the problem and the solution progress in inexact BD. \cite{mana2023surrogate} replaced the master problem with an RL-based surrogate to generate candidate solutions, without relying on repeatedly solving expensive master problems. 

In a different line of work, researchers have also investigated ML and RL approaches for general-purpose optimization algorithms. For example, ML has been leveraged to solve combinatorial optimization problems \citep{bengio2021machine,khalil2017learning,li2018combinatorial,khalil2016learning}. Recently, \cite{bello2017neural} proposed Neural Combinatorial Optimization, a framework that leverages policy gradient methods to train neural networks for generating high-quality solutions to traveling salesman problems. \cite{tang2020reinforcement} introduced a deep RL formulation to adaptively select cuts within a branch-and-cut framework for solving integer programs. \citet{chi2022deep} proposed a deep RL approach for column generation, where the column selection process is modeled as a sequential decision-making problem and optimized using Q-learning.

\section{Preliminaries}
In this section, we review the two-stage stochastic programming framework in Section \ref{sec:two-stage}, and the classical BD in Section \ref{sec:benders}, respectively.

\subsection{Two-Stage Stochastic Programming}\label{sec:two-stage}
We consider a class of two-stage stochastic programming problems, where in the first stage, a decision vector $x \in X \subseteq \mathbb{R}^{n_1}$ is chosen before the realization of uncertainty, and in the second stage, a recourse decision $y \in \mathbb{R}^{n_2}$ is made after observing the realization of the uncertain parameter $\xi$. The two-stage stochastic program can be written as $\min_{x \in X} c^{\mathsf T} x + \mathbb{E}_{\xi} [Q(x, \xi)]$,
where $c\in \mathbb{R}^{n_1}$ is the first-stage cost vector, and $\mathbb{E}_{\xi}[\cdot]$ denotes the expectation with respect to the distribution of $\xi$. Here, given the first-stage decision $x$ and uncertainty realization $\xi$, the second-stage recourse function $Q(x, \xi)$ is defined as
$Q(x, \xi) := \min_{y}\{q^{\mathsf T} y \ |\ Wy \ge h - Tx\}$,
where the uncertain parameter $\xi$ includes the data $(W, h, T, q)$.

We assume that the uncertain parameter $\xi$ can take values from a finite scenario set $\Omega$, where each scenario $\omega \in \Omega$ corresponds to a realization $\xi^\omega = (W^\omega, h^\omega, T^\omega, q^\omega)$ with probability $p^\omega$, satisfying $\sum_{\omega \in \Omega} p^\omega = 1$. Then, a sample average approximation (SAA) to the two-stage stochastic program can be written as 
\begin{equation}\label{eq:first-stage}
\min_{x \in X} \; c^{\mathsf T} x + \sum_{\omega \in \Omega} p^\omega Q^\omega(x),
\end{equation}
where
\begin{equation}\label{eq:second-stage}
Q^\omega(x) := \min_{y} \; (q^\omega)^{\mathsf T} y \quad 
\text{s.t.} \quad W^\omega y \ge h^\omega - T^\omega x.
\end{equation}

To ensure a well-defined problem, we further assume that the second-stage problem \eqref{eq:second-stage} is always feasible (i.e., satisfies relatively complete recourse) and that the recourse function is finite, i.e., $|Q^{\omega}(x)| < +\infty$ for all $x \in X$ and $\omega \in \Omega$. Our goal is to develop an efficient solution approach for the SAA problem \eqref{eq:first-stage}.


\subsection{Benders Decomposition}
\label{sec:benders}
The extensive-form SAA \eqref{eq:first-stage} becomes computationally challenging when the number of scenarios increases. BD exploits its separable structure by decomposing it into a master problem and a set of scenario-wise subproblems. 
At each iteration $t$, we first solve the following master problem $(\text{MP}_t)$:
\begin{align}
 (\text{MP}_t) \quad \min_{x \in X, \theta} \quad & c^{\mathsf T} x + \sum_{\omega \in \Omega} p^\omega \theta^\omega\nonumber\\
\text{s.t.} \quad &\theta^\omega \ge (\pi_{\tau}^{\omega})^{\mathsf T} (h^\omega - T^\omega x), \quad \forall \tau\in [t-1],\ \omega \in \Omega,\label{eq:master-cuts}
\end{align}
where $\theta^\omega \ge (\pi_{\tau}^{\omega})^{\mathsf T} (h^\omega - T^\omega x),\ \forall \tau\in [t-1],\ \omega \in \Omega$ denote the current optimality cuts to approximate the recourse function $Q^{\omega}(x)$ from below, and no cuts are present when $t=1$. The problem $(\text{MP}_t)$ is solved to obtain an optimal first-stage solution $(\hat{x}_t, \hat{\theta}_t)$ together with a lower bound to the original problem \eqref{eq:first-stage}. Given $\hat{x}_t$, we construct the following dual of the subproblem for each scenario $\omega\in\Omega$:
\begin{align}\label{eq:dual-sub}
 (\text{SP}_t^{\omega}) \quad &\max_{\pi^\omega\ge 0} \; (\pi^\omega)^{\mathsf T}(h^\omega - T^\omega \hat{x}_t)  \quad 
\text{s.t.} \quad (W^\omega)^{\mathsf T} \pi^\omega = q^\omega.
\end{align}
We denote the optimal solution of \eqref{eq:dual-sub} as $\pi_t^{\omega}$. By strong duality, $Q^{\omega}(\hat{x}_t)=(\pi_t^{\omega})^{\mathsf T}(h^\omega - T^\omega \hat{x}_t)$, and since $\hat{x}_t$ is a feasible first-stage solution, $c^{\mathsf T} \hat{x}_t+\sum_{\omega\in\Omega}p^{\omega}Q^{\omega}(\hat{x}_t)$ provides an upper bound for \eqref{eq:first-stage}. The optimal dual solutions $\{\pi_t^{\omega}\}_{\omega\in\Omega}$ are then used to generate Benders optimality cuts of the following form
\begin{equation}\label{eq:new-cut}
\theta^\omega \ge (\pi_t^{\omega})^{\mathsf T} (h^\omega - T^\omega x), \quad \forall \omega \in \Omega,
\end{equation}
which provide valid lower approximations of the recourse function. If the current $(\text{MP}_t)$'s optimal solution $(\hat{x}_t,\hat{\theta}_t)$ satisfies the newly generated cuts \eqref{eq:new-cut}, then we terminate the algorithm and return $(\hat{x}_t,\hat{\theta}_t)$ as the optimal solution. Otherwise, we add the optimality cuts \eqref{eq:new-cut} to the master problem and proceed to the next iteration. This iterative process converges to a globally optimal solution in a finite number of iterations under mild conditions, as shown in \cite{van1969shaped}.


Despite its effectiveness, classical BD often suffers from slow convergence. In particular, the number of generated cuts can grow rapidly, leading to computationally expensive master problems. Motivated by this computational challenge, we develop an RL framework to speed up the BD algorithm next.

\section{Reinforcement Learning for Benders Decomposition}
\begin{figure}[ht!]
    \centering
    \includegraphics[width=\linewidth]{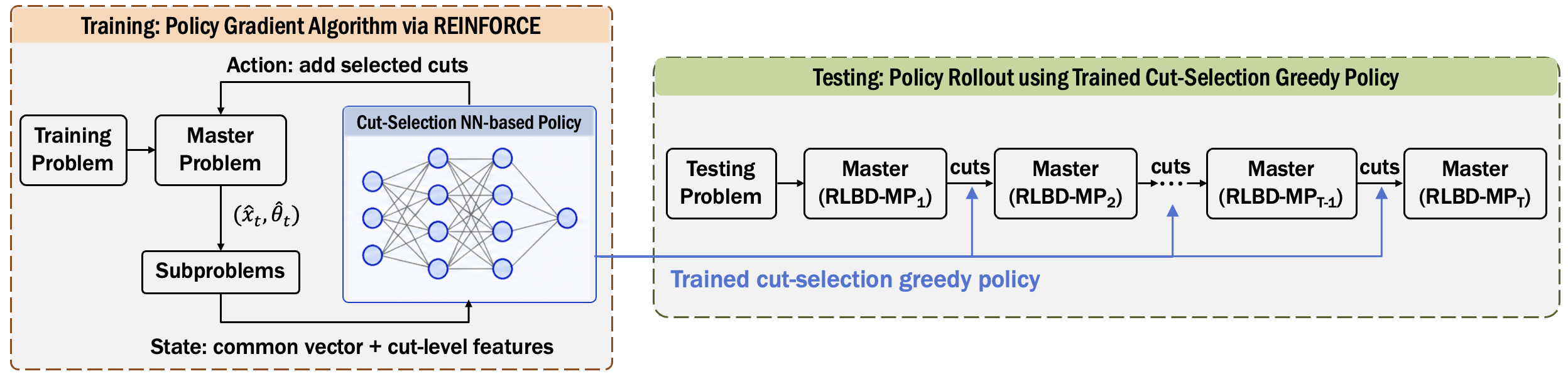}
    \caption{
    The RLBD Framework.
    }
    \label{fig:RLBD framework}
\end{figure}

In this section, we propose Reinforcement Learning for Benders Decomposition (RLBD), a framework that adaptively selects Benders cuts to improve computational efficiency. 
Unlike classical BD, which either adds all scenario-wise cuts or aggregates them into a single cut at each iteration, RLBD formulates cut selection as a sequential decision-making problem and learns a policy to select the most informative cuts at each iteration (see Figure \ref{fig:RLBD framework} for an illustration). 
We first introduce the MDP formulation in Section \ref{sec:MDP} and the neural network policy parameterization in Section \ref{sec:NN-policy}, respectively. The RLBD framework consists of two phases: (i) \textbf{Training}, where the stochastic cut-selection policy is learned using a policy gradient method via REINFORCE (Section \ref{sec:reinforce}), and (ii) \textbf{Testing}, where the trained policy is applied in a greedy manner to new problem instances via policy rollout (Section \ref{sec:testing}).

\vspace{-0.3cm}

\subsection{Markov Decision Process Formulation}
\label{sec:MDP}
We model the cut-selection problem in BD as a finite-horizon discounted MDP. In this formulation, each Benders iteration corresponds to one decision step. At iteration $t$, the agent observes the current master problem ($\text{RLBD-MP}_t$):
\begin{align}
 (\text{RLBD-MP}_t) \quad \min_{x \in X, \theta} \quad & c^{\mathsf T} x + \sum_{\omega \in \Omega} p^\omega \theta^\omega\nonumber\\
\text{s.t.} \quad &\theta^\omega \ge (\pi_{\tau}^{\omega})^{\mathsf T} (h^\omega - T^\omega x), \quad \forall \tau\in [t-1],\ \omega \in A_{\tau}.\label{eq:RLBD-master-cuts}
\end{align}
Unlike ($\text{MP}_t$) in classical BD \eqref{eq:master-cuts}, where all scenario-wise cuts $\omega\in\Omega$ are included, we instead maintain only a subset $A_{\tau}\subset \Omega$ for all $\tau\in[t-1]$. Moreover, the selected subset $A_{\tau}$ may vary across iterations. The agent solves ($\text{RLBD-MP}_t$) to obtain an optimal solution $(\hat{x}_t,\hat{\theta}_t)$ and an optimal objective value, which serves as a lower bound for problem \eqref{eq:first-stage}, denoted by $\text{LB}_t$. Based on $(\hat{x}_t,\hat{\theta}_t)$, the agent then solves the same subproblem ($\text{SP}_t^{\omega}$) in \eqref{eq:dual-sub} and generates scenario-wise cuts of the form \eqref{eq:new-cut} for all $\omega\in\Omega$. Since the first-stage solution $\hat{x}_t$ is feasible for \eqref{eq:first-stage} and $Q^{\omega}(\hat{x}_t)=(\pi_t^{\omega})^{\mathsf T}(h^\omega - T^\omega \hat{x}_t)$, the agent constructs an upper bound for problem \eqref{eq:first-stage}, defined by $\text{UB}_t:=c^{\mathsf T}\hat{x}_t+\sum_{\omega\in\Omega}p^{\omega}(\pi_t^{\omega})^{\mathsf T} (h^\omega - T^\omega \hat{x}_t)$.
Then, the agent selects a subset of candidate cuts $A_t\subset\Omega$, and the environment advances to the next iteration $t+1$ by adding the selected cuts $A_t$ to \eqref{eq:RLBD-master-cuts}.

Formally, the finite-horizon discounted MDP is defined by the tuple $(\mathcal{S}, \mathcal{A}, f, r,\gamma,T)$, where $\mathcal{S}$ is the state space, $\mathcal{A}$ is the action space, $f:\mathcal{S}\times\mathcal{A}\to\mathcal{S}$ represents the transition dynamics induced by the Benders solver, $r(s,a)$ is the step-wise reward function, $\gamma\in(0,1]$ is the discount factor, and $T$ is the number of decision stages. 
Next, we introduce each component in detail.

\subsubsection{State Space}
At each iteration $t$, the state representation contains both global information about the current convergence statistics, master problem and subproblem features, shared by all candidate cuts, and local information about each candidate cut. 

\textbf{Convergence statistics.}
    The state includes the current iteration index $t$, lower bound $\text{LB}_t$, upper bound $\text{UB}_t$, and the current optimality gap, defined by
    $\text{Gap}_t:=(\text{UB}_t-\text{LB}_t)/(|\text{UB}_t|+\epsilon)$. It also includes changes in the lower/upper bounds and the optimality gap, defined by
    $\Delta \text{LB}_t:=\text{LB}_t-\text{LB}_{t-1}$, $\Delta \text{UB}_t:=\text{UB}_{t-1}-\text{UB}_{t}$, and $\Delta \text{Gap}_t:=\text{Gap}_{t-1}-\text{Gap}_t$, respectively. Based on them, the state further includes the gap reduction rate $\rho_t^{\text{Gap}}:=\Delta \text{Gap}_t/(\text{Gap}_{t-1}+\epsilon)$, lower-bound improvement rate $\rho_t^{\text{LB}}:=\Delta \text{LB}_t/(|\text{LB}_{t-1}|+\epsilon)$, and upper-bound improvement rate $\rho_t^{\text{UB}}:=\Delta \text{UB}_t/(|\text{UB}_{t-1}|+\epsilon)$, which characterize the convergence statistics of the Benders procedure across iterations. Here, $\epsilon>0$ is a small positive constant, which is added in the denominator to avoid numerical instability. For each candidate cut of the form \eqref{eq:new-cut} generated at iteration $t$, we denote the cut-level violation at the current master solution $(\hat{x}_t,\hat{\theta}_t)$ as $v_t^{\omega}:=(\pi_t^{\omega})^{\mathsf T} (h^\omega - T^\omega \hat{x}_t)-\hat{\theta}_{t}^\omega$ for all $\omega\in\Omega$. We then record the mean violation $\bar v_t:=\sum_{\omega\in\Omega}p^{\omega}v_t^{\omega}$ and the maximum violation $v_t^{\max}:=\max_{\omega\in\Omega}v_t^{\omega}$, which reflect the approximation quality of the cuts. Note that $\bar v_t$ and $v
    _t^{\max}$ represent the global information, and we will include the cut-level violation $v_t^{\omega}$ in the local state component, introduced later.

\textbf{Master problem features.}
    The state includes the number of cuts added in the previous iteration $K$, the cumulative number of cuts added $C_t^{\mathrm{cum}}$, and the solve time for $(\text{RLBD-MP}_t)$, denoted by $T_t^{\text{MP}}$.
    These features reflect the size and computational burden of the master problem.

\textbf{Subproblem features.}
    After solving the scenario-wise subproblems \eqref{eq:dual-sub}, the state records summary statistics of the scenario-wise recourse costs.
    Let $Q_t^{\omega}$ denote the optimal objective value for ($\text{SP}_t^{\omega}$) in \eqref{eq:dual-sub} under the current first-stage solution $(\hat{x}_t,\hat{\theta}_t)$.
    We include the average recourse cost $\bar Q_t=\sum_{\omega\in\Omega}p^{\omega}Q_t^{\omega}$, the maximum recourse cost $Q_t^{\max}=\max_{\omega\in\Omega}Q_t^{\omega}$, the minimum recourse cost $Q_t^{\min}=\min_{\omega\in\Omega}Q_t^{\omega}$, and the standard deviation of the scenario-wise recourse costs $\{Q_t^{\omega}\}_{\omega\in\Omega}$.



\textbf{Cut-level features.}
    The cut-level state is given by
    $(v_t^{\omega},\|\pi_t^{\omega}\|_2,|(\pi_t^{\omega})^{\mathsf T} h^\omega|,\|(\pi_t^{\omega})^{\mathsf T} T^\omega\|_2,\text{NC}_t^{\omega})$ for each $\omega\in\Omega$,
    where $v_t^{\omega}$ is the cut-level violation, $\pi_t^{\omega}$ is the optimal dual variable used to generate the cut, $(\pi_t^{\omega})^{\mathsf T} h^\omega$ and $(\pi_t^{\omega})^{\mathsf T} T^\omega$ represent the intercept and coefficient of the cut, and
    $\text{NC}_t^{\omega}:=\sum_{\tau=1}^{t-1}\mathbf{1}\{\omega\in A_\tau\}$
    records the number of cuts generated by the same scenario $\omega$ in previous iterations, respectively.
    These local features capture the quality of each candidate cut, while the global features track the current convergence status.

The initial state $s_1$ is obtained from solving the master problem $(\text{RLBD-MP}_1)$ without any cuts. The initial $\text{LB}_1$, $\text{UB}_1$, $\text{Gap}_1$, cut violation statistics, and recourse-related statistics are computed from solving the initial master problem and subproblems. The history-dependent features, such as gap reduction and bound improvement rates, are initialized to zero because no previous iteration exists.

\subsubsection{Action Space, Transition Dynamics, and Reward Function}
\noindent\textbf{Action Space.}
At each iteration $t$, 
we use $\{1,2,\ldots,|\Omega|\}$ to encode all candidate cuts, each corresponding to one optimality cut of the form \eqref{eq:new-cut} with coefficients obtained from solving the subproblem \eqref{eq:dual-sub} for scenario $\omega\in\Omega$. The action corresponds to selecting a subset of $K$ cuts from $\{1,2,\ldots,|\Omega|\}$ to be added to the master problem. To this end, we model the action space as $\mathcal{A}=\{A\subset \{1,2,\ldots,|\Omega|\}:\ |A|=K\}$. Note that the same cut index $a_t\in \Omega$ may correspond to different cuts in different iterations $t$, because even for the same scenario, the cut coefficients obtained from solving the subproblem \eqref{eq:dual-sub} can vary depending on the first-stage solution.

We use $\pi_{\theta}(a_t|s_t)$ to denote the probability of selecting cut $a_t\in\{1,2,\ldots,|\Omega|\}$ given state $s_t$. Note that since the iteration index $t$ is included in state $s_t$, this policy $\pi_{\theta}(a_t|s_t)$ is naturally time-dependent.
Given state $s_t$, we sample an ordered sequence $A_t=(a_t^{(1)},a_t^{(2)},\ldots,a_t^{(K)})$ following $\pi_{\theta}(\cdot|s_t)$ without replacement. The sampling is performed sequentially: after one cut is selected, it is masked out, and the remaining probabilities are re-normalized before the next cut is sampled. The joint probability of action $A_t\in\mathcal{A}$ is thus given by $P_{\theta}(A_t|s_t)=\prod_{i=1}^K \frac{\pi_{\theta}(a_t^{(i)}|s_t)}{1-\sum_{j<i}\pi_{\theta}(a_t^{(j)}|s_t)}$.

This stochastic cut-selection policy naturally introduces exploration through the policy distribution. Unlike deterministic selection rules, the use of a stochastic policy allows the agent to sample actions according to their relative importance, assigning higher probability to more promising cuts while retaining the ability to explore alternative candidates. We will introduce a neural network architecture to parameterize $\pi_{\theta}(a_t|s_t)$ based on both global information and cut-level features in Section \ref{sec:NN-policy}.

\noindent\textbf{Transition Dynamics.}
Given the current state $s_t$ and action $A_t$, the environment advances to the next state $s_{t+1}=f(s_t,A_t)$ by adding the selected cuts $A_t$ to the master problem \eqref{eq:RLBD-master-cuts}.
The updated master problem $(\text{RLBD-MP}_{t+1})$ is re-optimized to obtain a new first-stage solution. Based on this solution, we solve the subproblems $(\text{SP}^{\omega}_{t+1})$ for all scenarios $\omega\in\Omega$ and generate a new set of candidate cuts. The resulting information from the master problem and subproblems, including updated bounds, optimality gap, violation, and other summary statistics, is used to construct the next state $s_{t+1}$.


\noindent\textbf{Reward Function.}
The step-wise reward is defined as $r(s_t,A_t) := \alpha \cdot F_t + \beta \cdot L_t - \lambda$,
where $F_t$ denotes the gap-based reward term, $L_t$ denotes the reward associated with the master problem solve time, and $-\lambda$ is a per-step penalty. 
First of all, the gap-based reward is defined by $F_t := \log(\text{Gap}_{t-1}) - \log(\text{Gap}_t)$, which is scale-invariant and designed to reflect the convergence progress of BD.
Second, $L_t := - \frac{T^{\mathrm{MP}}_{t}}{T_{\mathrm{ref}}}$,
where $T^{\mathrm{MP}}_{t}$ is the master problem solve time at iteration $t$, and $T_{\mathrm{ref}}$ is a normalization constant chosen based on the problem scale. 
Finally, the per-step penalty $-\lambda$ is introduced to encourage convergence within fewer iterations, which accumulates over time and imposes an increasing penalty on prolonged solution trajectories. Note that the parameters $\alpha$, $\beta$ and $\lambda$ are user-defined weights that balance these components and will be tuned in Section \ref{sec:tuning}.




Combining all components above, the value function of this MDP takes the following form: 
$$\max_{\theta} V_{\theta}(s_1)
:=\mathbb{E}\left[\sum_{t=1}^T\gamma^tr(s_t,A_t)\right],\ \text{where } A_t\sim P_{\theta}(\cdot|s_t),\ s_{t+1}=f(s_t,A_t),\ \forall t\in [T].$$

\subsection{Neural Network Policy Parameterization}
\label{sec:NN-policy}
\begin{figure}[ht!]
    \centering
    \includegraphics[width=\linewidth]{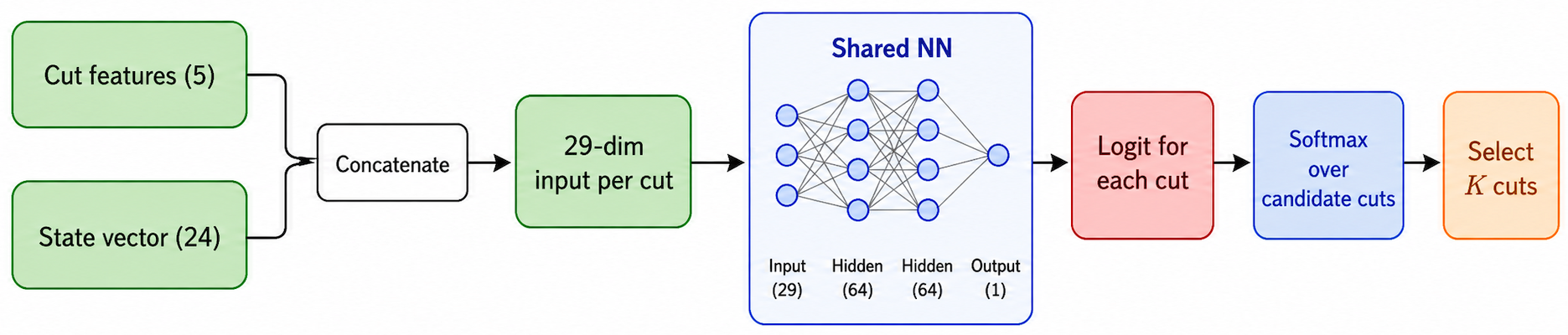}
    \caption{
    The shared neural network is applied to each candidate cut using the concatenation of common state vector and local cut features, producing one logit per candidate cut.
    }
    \label{fig:NN_Parameterization}
\end{figure}

We parameterize the cut-selection policy $\pi_{\theta}(a_t|s_t)$ using a shared feedforward neural network (see Figure \ref{fig:NN_Parameterization}). The purpose of the network is to assign a scalar score to each candidate cut, which is later converted into a probability distribution for stochastic cut selection. Specifically, at iteration $t$, for each candidate cut $\omega \in \Omega$, we construct a cut-specific state representation $s_t^{\omega}\in\mathbb{R}^{28}$. The first 23 entries of $s_t^{\omega}$ include the global state information shared by all candidate cuts, and the last 5 entries describe the local properties of candidate cut $\omega$, given by $(v_t^{\omega},\|\pi_t^{\omega}\|_2,|(\pi_t^{\omega})^{\mathsf T} h^\omega|,\|(\pi_t^{\omega})^{\mathsf T} T^\omega\|_2,\text{NC}_t^{\omega})$.

Each input vector $s_t^{\omega}$ is passed through a multilayer perceptron with two hidden layers and ReLU activations: $z_\omega= \mathrm{NN}_{\theta}(s_t^{\omega})=W_3\cdot\mathrm{ReLU}(W_2\cdot\mathrm{ReLU}(W_1s_t^{\omega}+b_1)+b_2)+b_3$, for all $\omega\in\Omega$,
where $\theta=(W_l,b_l)_{l=1}^3$ denotes the network parameters, and $z_\omega \in \mathbb{R}$ is the output logit assigned to candidate cut $\omega$. Given the logits produced by the neural network, we apply a softmax function to transform them into a probability distribution: $\pi_\theta(\cdot|s_t)=\text{Softmax}(z_1,z_2,\ldots,z_{|\Omega|})$.
Here, $\pi_{\theta}(\omega|s_t)$ denotes the probability of selecting cut $\omega$, which is used to construct the joint action probability $P_{\theta}(A_t|s_t)$.

\subsection{Training: Policy Gradient Algorithm via REINFORCE}
\label{sec:reinforce}
In this section, we leverage policy gradient methods to train the stochastic cut-selection policy $\pi_{\theta}$ via the REINFORCE algorithm \citep{Williams1992} (see Algorithm \ref{alg:reinforce_benders}). One training episode corresponds to one complete BD trajectory. During each episode, at each iteration $t$, the algorithm records the current state $s_t$, the selected cuts $A_t$, their log-probabilities $\log P_{\theta}(A_t|s_t)$, and the stepwise reward $r_t$. After the episode terminates (either BD converges or hits the maximum iterations $T_{\max}$), the return-to-go is computed for each step: $G_t =
\sum_{\ell=t}^{T}
\gamma^{\ell-t} r_\ell$. According to \cite{Williams1992}, the gradient of the value function $V_{\theta}(s_1)$ can be estimated as $\nabla_{\theta}V_{\theta}(s_1)\approx\mathbb{E}[\sum_{t=1}^T\nabla_{\theta}\log P_{\theta}(A_t|s_t)\cdot G_t]$,
where $\nabla_{\theta}\log P_{\theta}(A_t|s_t)=\sum_{i=1}^K\nabla_{\theta}\log(\frac{\pi_{\theta}(a_t^{(i)}|s_t)}{1-\sum_{j<i}\pi_{\theta}(a_t^{(j)}|s_t)})$. At the end of each episode, we apply backpropagation to update the policy network parameter $\theta$ in the gradient ascent direction. 





\begin{algorithm}[t]
\small
\caption{Training Cut-Selection Policy via REINFORCE}
\label{alg:reinforce_benders}

\SetKwInOut{Input}{Input}
\SetKwInOut{Output}{Output}

\Input{Training problem, initial policy network parameter $\theta$, discount factor $\gamma$, number of cuts added per iteration $K$, maximum iterations $T_{\max}$, number of episodes $E$, learning rate $\eta$, optimality gap tolerance $\epsilon_{\text{tol}}$, reward weights $\alpha,\ \beta,\ \lambda$, reference master problem solve time $T_{\text{ref}}$}
\Output{Trained cut-selection policy $\pi_\theta$}


\For{$e = 1,\ldots,E$}{
    Initialize the master problem $(\text{RLBD-MP}_1)$ based on the training problem without any cuts\;
    
    \For{$t = 1,\ldots,T_{\max}$}{
        Solve $(\text{RLBD-MP}_t)$ to obtain an optimal solution $(\hat{x}_t,\hat{\theta}_t)$ and a lower bound $\text{LB}_t$\;

        Based on $(\hat{x}_t,\hat{\theta}_t)$, solve $(\text{SP}^{\omega}_{t})$ to construct candidate cuts \eqref{eq:new-cut} and obtain an upper bound $\text{UB}_t$\;
        
        Construct state $s_t$ according to the global information and cut-level features\;
        
        
        Compute policy distribution
        $\pi_\theta(\cdot|s_t)=\text{Softmax}(\text{NN}_\theta(s_t^1),\ldots,\text{NN}_\theta(s_t^{|\Omega|}))$\;
        Sample $K$ cuts following $\pi_\theta(\cdot|s_t)$ without replacement and set the selected subset as action $A_t$\;

        Compute the stepwise reward $r_t=r(s_t,A_t)= \alpha \cdot F_t + \beta \cdot L_t - \lambda$\;
        Add the selected cuts $A_t$ to the master problem and advance to the next iteration\;
        
        \If{$(\text{UB}_t-\text{LB}_t)/|\text{UB}_t|<\epsilon_{\text{tol}}$}{
            break\;
        }
    }
    
    Compute the return-to-go backward: $G_T=0$ and
    $G_t = r_t + \gamma G_{t+1}$ for $t=T-1,\ldots,0$\;
    
    
    Update the policy parameter:
    $\theta \leftarrow \theta+\eta\sum_{t=1}^{T}
    \nabla_\theta  \log P_\theta(A_t\mid s_t)\cdot G_t$\;
}

\Return{$\pi_\theta$.}

\end{algorithm}

\subsection{Testing: Policy Rollout using the Trained Cut-Selection Greedy Policy}\label{sec:testing}
Once the cut-selection policy $\pi_{\theta}$ is trained, it can be applied to solve problems with similar structures but different inputs. We consider a setting in which a decision maker needs to repeatedly solve a collection of two-stage stochastic programs that share the same problem structure but differ in their realizations of uncertain parameters and/or dimensions of the decision variables. This setting naturally arises in applications where similar decision problems must be solved under different scenarios or operating conditions. For example, in supply chain management, one may solve a sequence of stochastic facility location problems that share the same objective function and constraints but differ in demand realizations or  decision dimensions. Since the policy network $\pi_{\theta}$ is designed to be independent of the problem size, it can naturally generalize to such instances with different demand realizations and varying numbers of candidate facilities and customer sites. For a given test problem, we then implement BD with the trained cut-selection policy in a greedy manner: at each iteration, we select the top-$K$ cuts with the largest logits from the policy network to add to the master problem.





    
    
    



\section{Numerical Results}
We test our RLBD on a two-stage stochastic electric vehicle (EV) charging station location problem and compare it with existing benchmarks. We first introduce our experimental setup in Section \ref{sec:setup}, report the hyperparameter tuning results in Section \ref{sec:tuning}, and compare our RLBD with other benchmarks in Section \ref{sec:benchmark}, respectively. Finally, we provide some managerial insights from the trained cut-selection policy in Section \ref{sec:insight}.
\subsection{Experimental Setup}\label{sec:setup}
We consider a two-stage stochastic EV charging station location problem, where we have a set $\mathcal I$ of candidate charging stations and a set $\mathcal J$ of customer sites with uncertain demand $d_j,\ j\in\mathcal{J}$. In the first stage, a subset of charging stations is selected to open, and their capacities are determined before demand realization. In the second stage, demand is realized and fulfilled by the open stations. The objective is to minimize first-stage investment costs and second-stage transportation and unmet-demand penalties, while maximizing total revenue. 
In the training phase, we generate $|\mathcal I|=8$ candidate charging stations and $|\mathcal J|=12$ customer sites. For each customer site $j\in\mathcal{J}$, the demand $d_j$ follows a Normal distribution. We first sample $\mu_j$ uniformly from $[20,100]$ and set $\sigma_j=0.1\mu_j$ and then generate $|\Omega|=100$ demand scenarios $\{d_j^{\omega}\}_{\omega\in\Omega}$ following $\mathcal{N}(\mu_j,\sigma_j^2)$. This forms our training problem. Please refer to \eqref{eq:obj}--\eqref{eq:cons5} in Appendix \ref{append:EV} for the full model formulation and parameter setup.

We train the policy network using the \texttt{Adam} optimizer, with learning rate $10^{-3}$, and discount factor $\gamma=0.99$. All experiments are conducted on the Pitzer cluster at the Ohio Supercomputer Center with 16 CPU cores. We use Gurobi 12.0.0 coded in Python 3.10 for solving all optimization problems.

\subsection{Training: Hyperparameter Tuning}\label{sec:tuning}
We first tune the hyperparameters in defining the reward function. We fix $\beta = 0.001$ and $T_{\mathrm{ref}} = 0.1$ and perform a grid search over the other weight parameters $\alpha\in \{0.01, 0.1, 1.0\}$ and $\lambda\in \{0.001, 0.01, 0.1\}$, resulting in nine reward configurations. For each configuration, we train the policy network in five independent runs and evaluate the learned cut-selection policies in a greedy manner on the same training problem. The evaluation results are summarized in Figure~\ref{fig:reward_hyperparameter_tuning}, where we plot the episodic elapsed time, final optimality gap, and number of Benders iterations across different reward configurations. The shaded area represents the 90\% confidence interval among the five independent runs. Based on these curves, $\alpha = 0.01$ and $\lambda = 0.001$ perform the best in all three metrics, and we fix $\alpha = 0.01$ and $\lambda = 0.001$ in the following experiments.

After fixing the reward setting, we now tune the number of cuts added per iteration $K$ from $\{10,20,30\}$. 
The average evaluation results over five independent runs are shown in Figure~\ref{fig:k_sensitivity}, where we report the overall trends of  elapsed time, final optimality gap, and the number of Benders iterations under different values of $K$. Note that in these experiments, even under the same initialized policy, different values of $K$ lead to different BD performance, and therefore, the three curves do not start from the same initial point. We also provide additional results of total reward and their close-up in Figures~\ref{fig:k_tuning_results_append} in Appendix \ref{append:results}. Overall, $K=10$ achieves the most stable convergence behavior and computational efficiency, and we adopt $K=10$ in the next benchmark experiments.

\newcommand{\rewardlegendwidth}{0.45\textwidth}

\begin{figure*}[htbp]
    \centering
\includegraphics[width=\rewardlegendwidth]{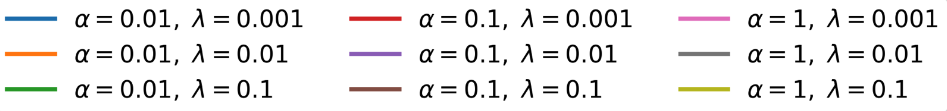}

    \begin{subfigure}[b]{0.32\textwidth}
        \centering
        \includegraphics[width=\textwidth]{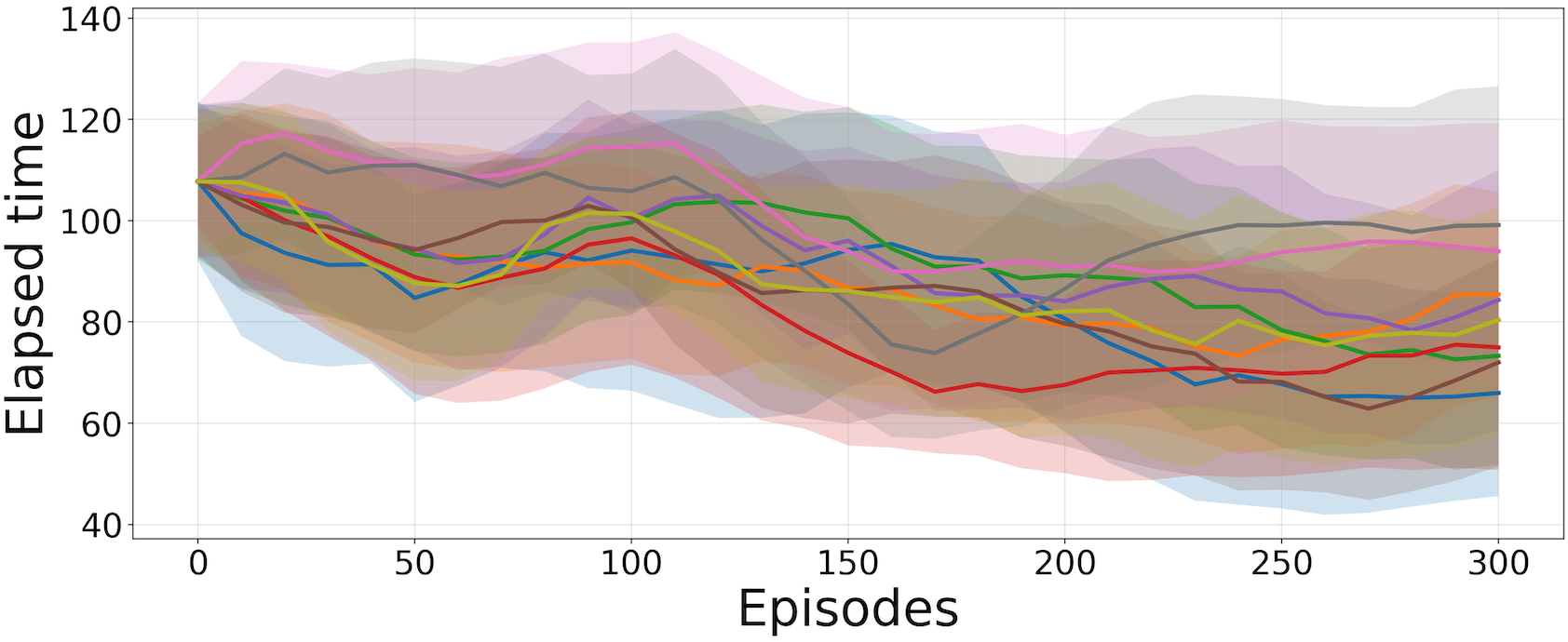}
        \caption{Elapsed Time}
        \label{fig:reward_time}
    \end{subfigure}
    \hfill
    \begin{subfigure}[b]{0.32\textwidth}
        \centering
        \includegraphics[width=\textwidth]{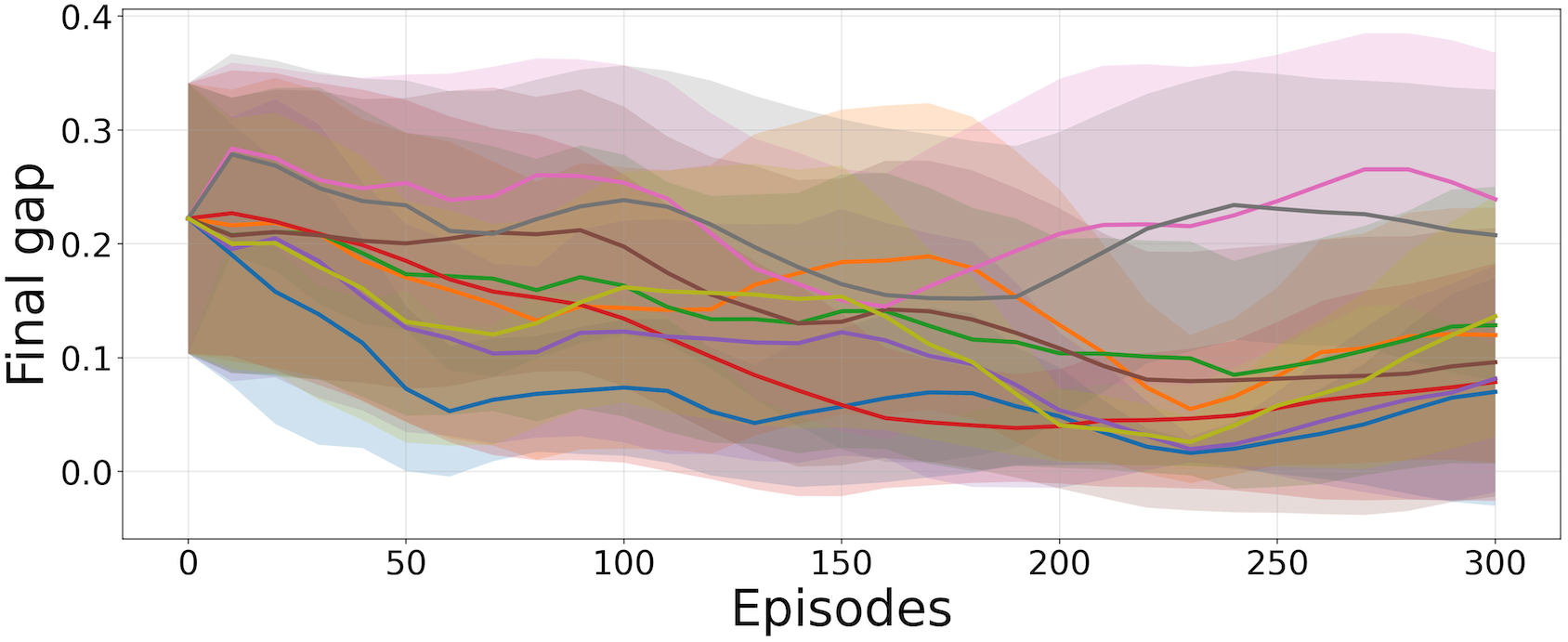}
        \caption{Final Optimality Gap}
        \label{fig:reward_gap}
    \end{subfigure}
    \hfill
    \begin{subfigure}[b]{0.32\textwidth}
        \centering
        \includegraphics[width=\textwidth]{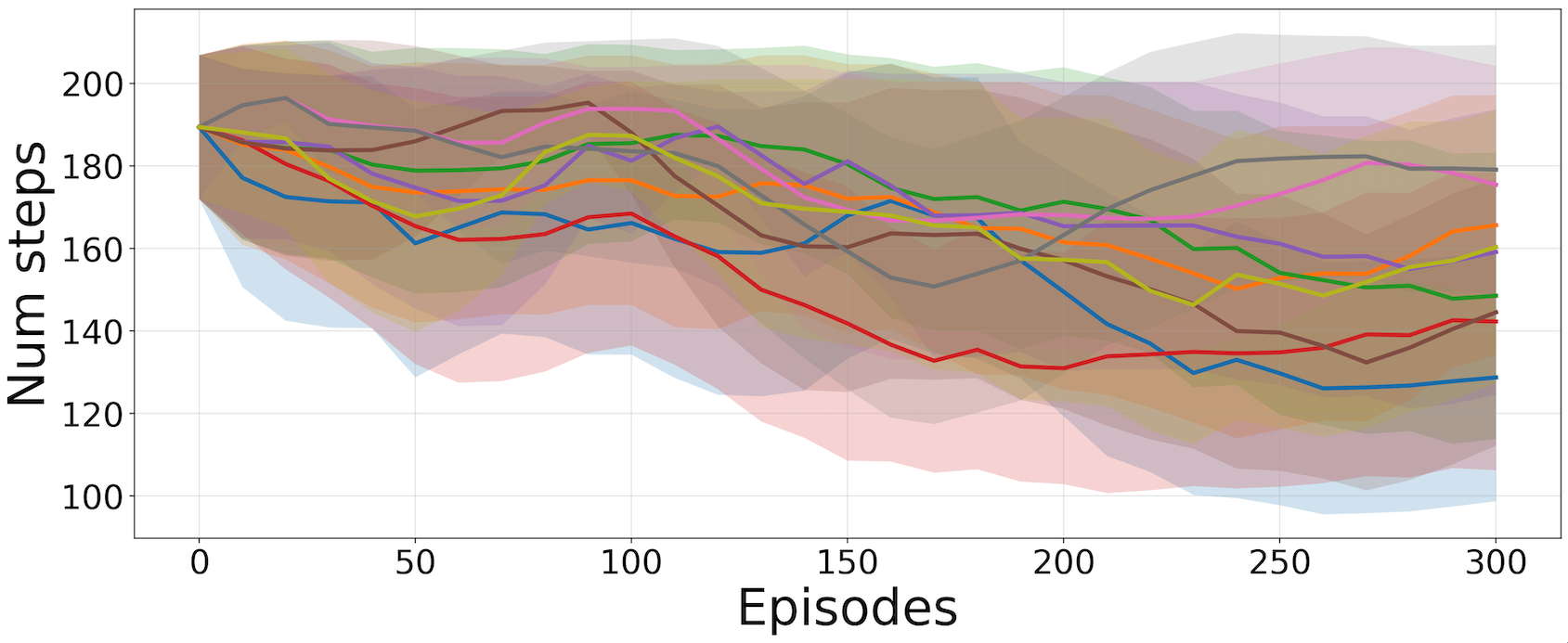}
        \caption{Number of Steps}
        \label{fig:reward_steps}
    \end{subfigure}


    \caption{Hyperparameter tuning results under different combinations of $\alpha$ and $\lambda$.}
    \label{fig:reward_hyperparameter_tuning}
\end{figure*}

\newcommand{\klegendwidth}{0.25\textwidth}

\begin{figure*}[htbp]
    \centering
\includegraphics[width=\klegendwidth]{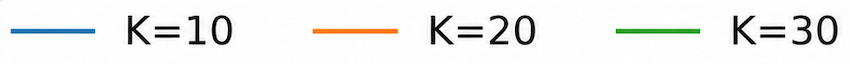}

    \begin{subfigure}[b]{0.32\textwidth}
        \centering
        \includegraphics[width=\textwidth]{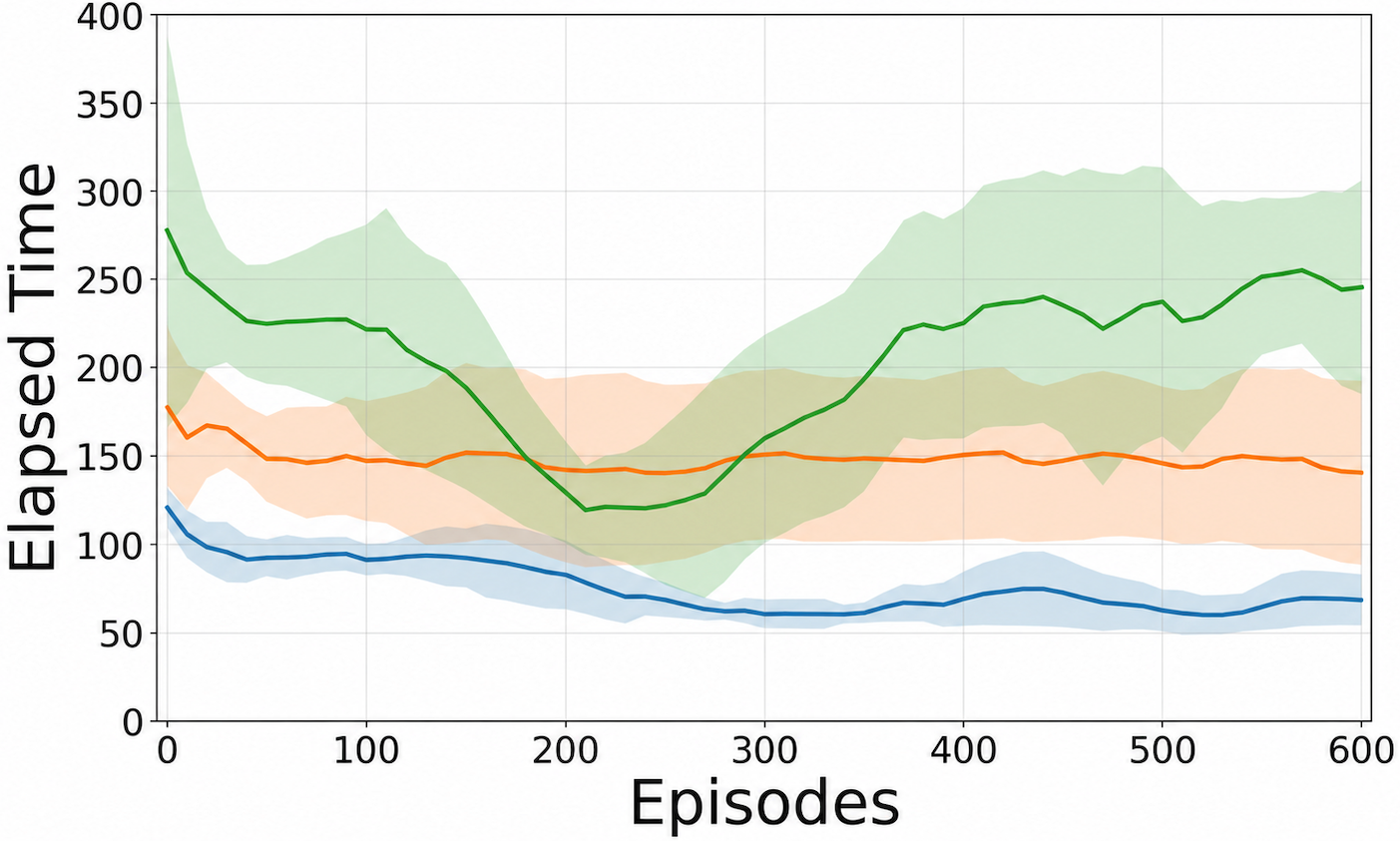}
        \caption{Elapsed Time}
        \label{fig:k_time}
    \end{subfigure}
    \hfill
    \begin{subfigure}[b]{0.32\textwidth}
        \centering
        \includegraphics[width=\textwidth]{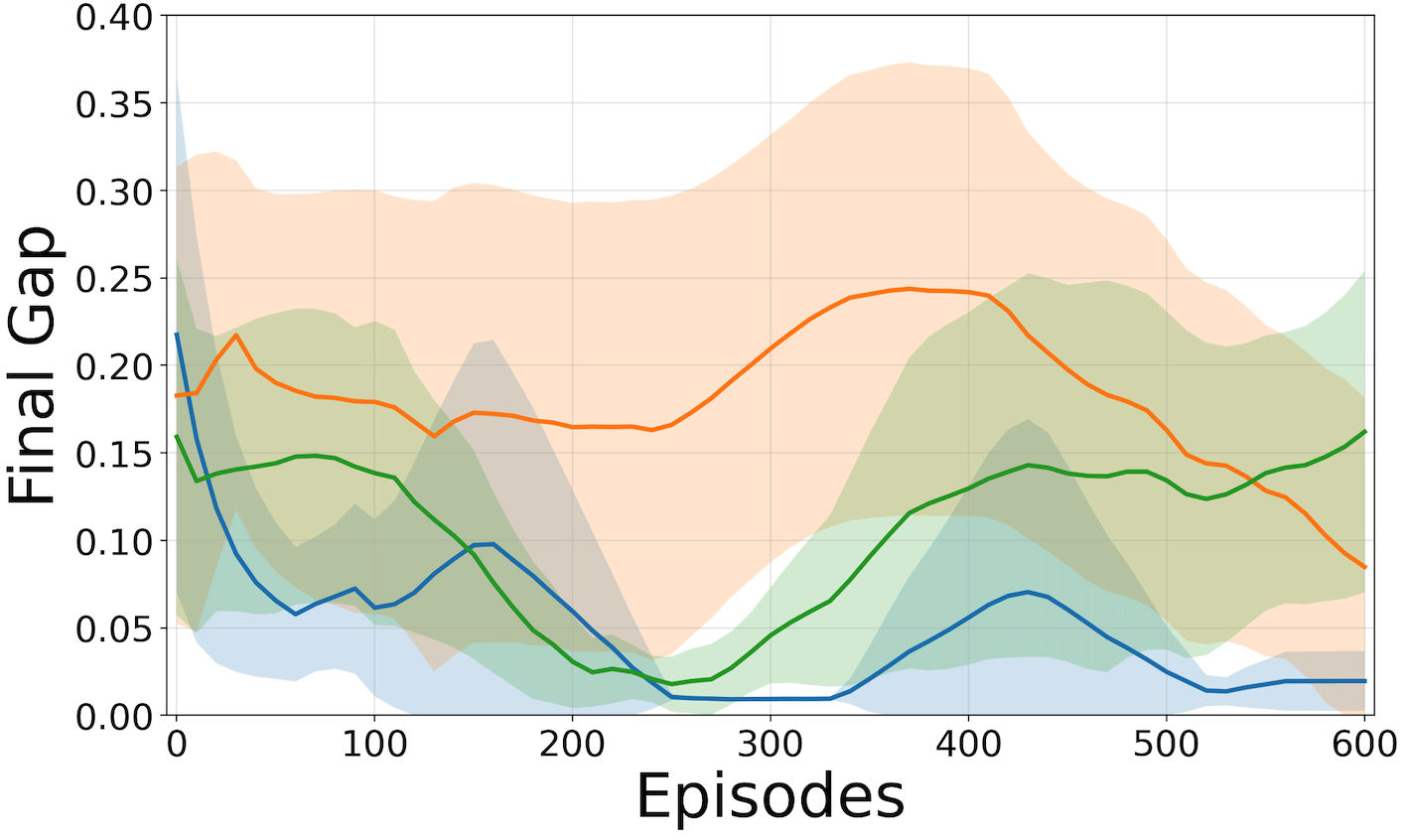}
        \caption{Final Optimality Gap}
        \label{fig:k_gap}
    \end{subfigure}
    \hfill
    \begin{subfigure}[b]{0.32\textwidth}
        \centering
        \includegraphics[width=\textwidth]{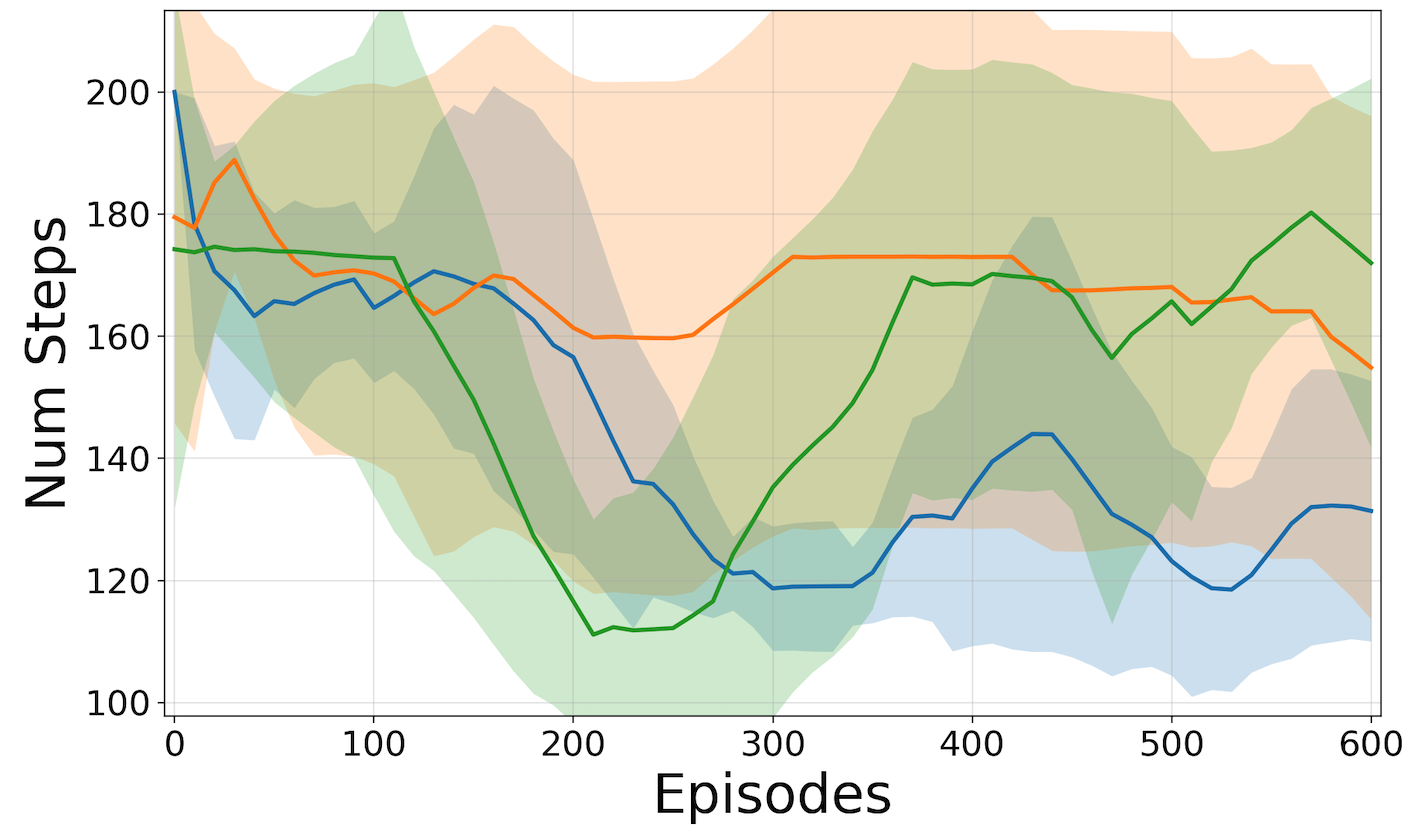}
        \caption{Number of Steps}
        \label{fig:k_steps}
    \end{subfigure}


    \caption{Evaluation results under different numbers of cuts added per iteration $K$.}
    \label{fig:k_sensitivity}
\end{figure*}

\subsection{Testing: Comparison with Existing Benchmarks}\label{sec:benchmark}
Once the cut-selection policy $\pi_{\theta}$ is trained in our RLBD framework, we apply it to solve testing instances with various demand distributions and problem sizes and compare it with the following benchmarks: (i) vanilla BD presented in Section \ref{sec:benders}; and (ii) LearnBD proposed in \cite{jia2021benders}, where we train a support vector machine classifier on the same training problem. 
We consider several distribution-shift settings in which we sample $|\Omega|=100$ testing demand scenarios from the same Normal distribution as the training problem, as well as from left-skewed and right-skewed Normal distributions.

Table~\ref{tab:benchmark_results} summarizes the comparison between BD, LearnBD, 
and RLBD on different testing instances, where we report the average solution time (Time), master problem solve time (Master), number of iterations (Iter.), and the final optimality gap (Gap) over five independent testing datasets. The time limit is set to one hour and the optimality gap tolerance is set to $\epsilon_{\text{tol}}=1\%$ for all instances. From Table~\ref{tab:benchmark_results}, when $|\mathcal I|\times |\mathcal J|=8\times 12$ or $10\times 15$, all methods solve the problem to the prescribed optimality gap tolerance, while RLBD achieves approximately a fivefold reduction in computational time. This improvement is primarily due to smaller master problems with adaptive cut selection (see the ``Master'' column). When the problem size increases to $|\mathcal I|\times |\mathcal J|=20\times 30$, none of the methods can solve the problem within one hour; however, RLBD attains a smaller final optimality gap.

\begin{table*}[htbp]
\centering
\caption{Performance comparison between BD, LearnBD, and RLBD on testing instances with various demand distributions and problem sizes}
\label{tab:benchmark_results}
\resizebox{\textwidth}{!}{
\begin{tabular}{llrrrrrrrrrrrr}
\toprule
\multirow{2}{*}{\makecell{$|\mathcal I|\times |\mathcal J|$}}
& \multirow{2}{*}{Method}
& \multicolumn{4}{c}{Normal}
& \multicolumn{4}{c}{Left-skewed Normal}
& \multicolumn{4}{c}{Right-skewed Normal} \\
\cmidrule(lr){3-6} \cmidrule(lr){7-10} \cmidrule(lr){11-14}
&
& Time (s) & Master (s) & Iter. & Gap (\%)
& Time (s) & Master (s) & Iter. & Gap (\%)
& Time (s) & Master (s) & Iter. & Gap (\%) \\
\midrule

\multirow{4}{*}{$8 \times 12$}
& BD
& 237.63 & 222.03 & 90.2 & 0.55
& 244.52 & 227.63 & 90.6 & 0.57
& 252.60 & 234.63 & 90.4 & 0.56 \\

& LearnBD
& 234.77 & 212.56 & 124.2 & 0.55
& 239.75 & 215.14 & 125.0 & 0.57
& 246.73 & 221.81 & 124.6 & 0.56 \\


& \textbf{RLBD}
& \textbf{44.43} & \textbf{25.62} & 106.6 & 0.95
& \textbf{46.00} & \textbf{26.39} & 108.0 & 0.90
& \textbf{48.32} & \textbf{27.16} & 106.6 & 0.87 \\

\midrule

\multirow{3}{*}{$10 \times 15$}
& BD
& 446.70 & 418.17 & 121.4 & 0.69
& 436.67 & 404.55 & 118.4 & 0.87
& 454.84 & 423.42 & 122.2 & 0.82 \\

& LearnBD
& 448.74 & 408.87 & 163.8 & 0.71
& 440.99 & 396.62 & 161.2 & 0.87
& 453.92 & 410.03 & 164.2 & 0.82 \\


& \textbf{RLBD}
& \textbf{97.37} & \textbf{60.25} & 154.4 & 0.86
& \textbf{102.26} & \textbf{62.00} & 151.4 & 0.97
& \textbf{102.77} & \textbf{62.66} & 150.8 & 0.79 \\

\midrule

\multirow{3}{*}{$20 \times 30$}
& BD
& 3600 & 3600 & 185.0 & 38.89
& 3600 & 3600 & 185.8 & 35.48
& 3600 & 3600 & 186.4 & 39.07 \\

& LearnBD
& 3600 & 3600 & 210.2 & 38.64
& 3600 & 3600 & 210.2 & 35.45
& 3600 & 3600 & 211.6 & 39.07 \\


& \textbf{RLBD}
& 3600 & 3600 & 325.2 & \textbf{28.39}
& 3600 & 3600 & 329.0 & \textbf{29.41}
& 3600 & 3600 & 332.0 & \textbf{29.21} \\

\bottomrule
\end{tabular}
}
\end{table*}

\subsection{Managerial Insights from the Learned Cut-Selection Policy}\label{sec:insight}
\begin{figure*}[htbp]    \centering    \begin{subfigure}[b]{0.32\textwidth}        \centering        \includegraphics[width=\textwidth]{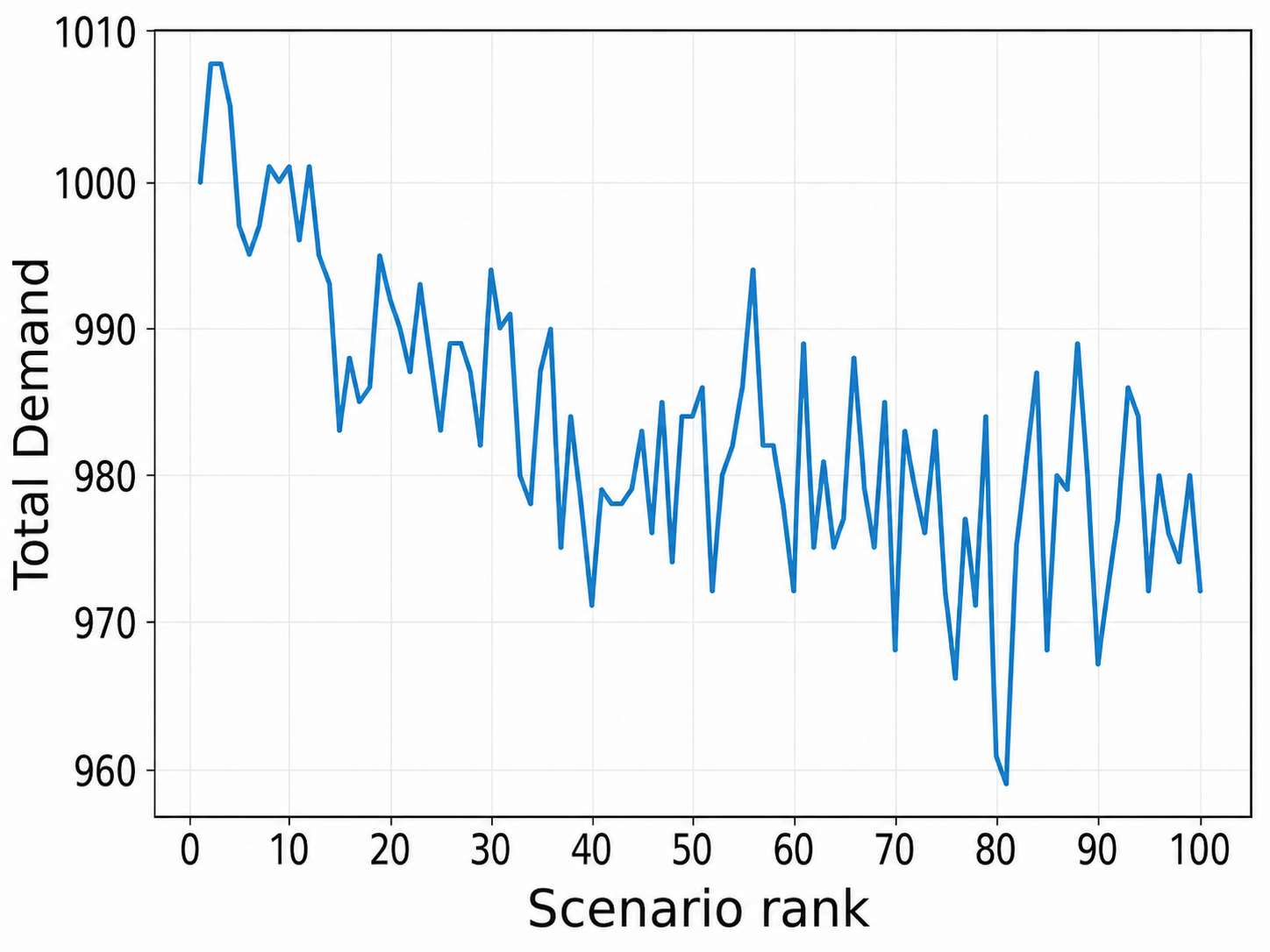}        \caption{Total Demand }        \label{fig:scenario_total_demand}    \end{subfigure}    \hfill    \begin{subfigure}[b]{0.32\textwidth}        \centering        \includegraphics[width=\textwidth]{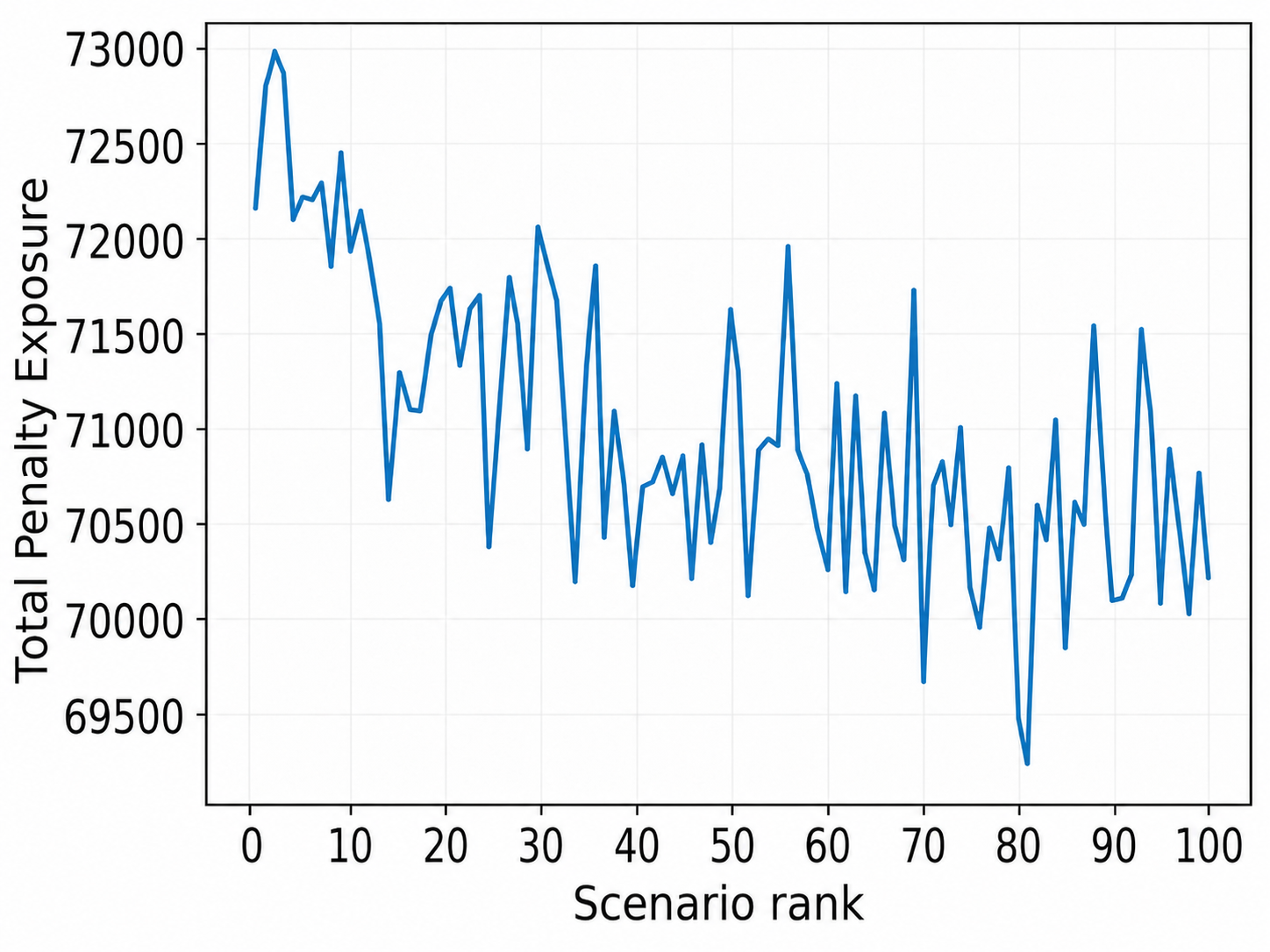}        \caption{Total Penalty Exposure}        \label{fig:scenario_total_penalty_exposure}    \end{subfigure}    \hfill    \begin{subfigure}[b]{0.32\textwidth}        \centering        \includegraphics[width=\textwidth]{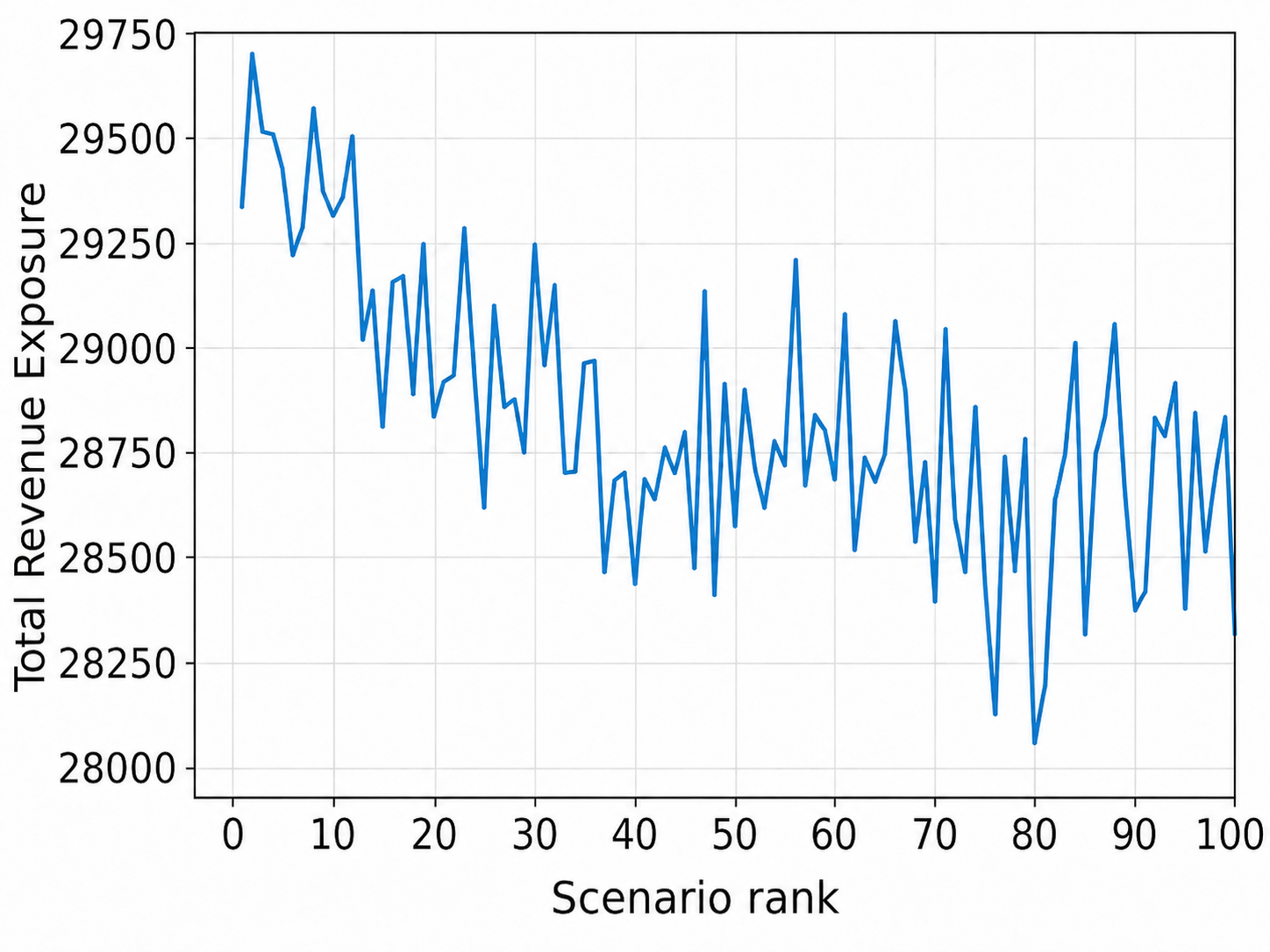}        \caption{Total Revenue Exposure}        \label{fig:scenario_total_revenue_exposure}    \end{subfigure}    \caption{Scenario-level exposure metrics sorted by scenario rank.}    \label{fig:scenario_exposure_metrics}\end{figure*}

To understand the relationship between the behavior of the trained cut-selection policy $\pi_{\theta}$ in RLBD and the characteristics of the input dataset, we rank the scenarios based on the number of cuts generated by each scenario across all iterations ($\text{NC}_T^{\omega}$) and plot their corresponding total demand $\sum_{j\in\mathcal J}d_j^{\omega}$, total penalty exposure $\sum_{j\in\mathcal J}p_jd_j^{\omega}$, and total revenue exposure $\sum_{j\in\mathcal J}r_jd_j^{\omega}$ in Figure \ref{fig:scenario_exposure_metrics}. On the $x$-axis, scenarios are sorted by the selected count $\text{NC}_T^{\omega}$ in descending order, where Scenario \#1 corresponds to the most frequently selected scenario. We focus on the testing instance with right-skewed Normal distribution and $|\mathcal I|\times|\mathcal J|=10\times 15$. From Figure \ref{fig:scenario_exposure_metrics}, the trained policy tends to select cuts from scenarios with higher total demand, penalty exposure, and revenue exposure, suggesting that high-load demand realizations correspond to more critical scenarios and are more likely to generate informative cuts.


\section{Conclusions}
In this work, we proposed an RLBD framework for adaptively selecting Benders cuts to improve the efficiency of solving two-stage stochastic programs. The cut-selection process was modeled as an MDP, and the stochastic cut-selection policy was trained using the REINFORCE algorithm. The learned policy can be deployed on testing instances with varying uncertainty distributions and problem sizes. Numerical experiments on a stochastic EV charging station location problem demonstrated that RLBD generalizes effectively to problems with similar structures but different data inputs (including distribution shifts) and decision variable dimensions. Compared with existing benchmark methods, RLBD achieves up to a fivefold reduction in computational time in small to medium testing instances and obtains smaller optimality gaps in larger testing instances.

\bibliographystyle{plainnat}
\bibliography{reference}

\newpage
\appendix
\begin{center}
    \Large \textbf{Appendix}
\end{center}

\section{Two-stage Stochastic EV Charging Station Location problem}\label{append:EV}
Let $\mathcal I$, $\mathcal J$, and $\Omega$ denote the set of candidate charging stations, customer sites, and demand scenarios, respectively. 
Define the first-stage decision variables as $y_i\in\{0,1\}$ and $z_i\in\mathbb Z_+$ for all $i\in\mathcal I$, where $y_i=1$ if station $i$ is selected to open and $z_i$ denotes the number of chargers installed at station $i$. Define the second-stage decision variables $x_{ij}$ and $s_j$ as the amount of demand at customer site $j \in \mathcal J$ fulfilled by charging station $i$ and the unmet demand at customer site $j$ for all $i \in \mathcal I, \ j\in \mathcal J$, respectively. Parameters $f_i,b_i, c_{ij},p_j,r_j, C_i$ denote the investment cost for opening station $i\in I$, unit cost for installing one charger at station $i\in I$, unit transportation cost from station $i\in I$ to customer site $j\in J$, unit penalty for unmet demand at customer site $j\in J$, unit revenue for demand fulfillment at site $j\in J$, and the maximum number of EV demand each charger can fulfill at location $i\in I$, respectively.

Following \cite{sun2025contextual}, the two-stage stochastic EV charging station location model is given by
\begin{align}
\min_{y,z}
&\quad
\sum_{i\in\mathcal I} f_i y_i
+
\sum_{i\in\mathcal I} b_i z_i
+
\frac{1}{|\Omega|}
\sum_{\omega\in\Omega} Q^\omega(z)\label{eq:obj}\\
\text{s.t.}& \quad z_i \le M_i y_i,\ \forall i\in\mathcal I \label{eq:cons1}\\
& \quad y_i\in\{0,1\},\ z_i\in\mathbb Z_+,\ \forall i\in\mathcal I \label{eq:cons2}
\end{align}
where constraints \eqref{eq:cons1} enforce that chargers can be installed at location \( i\in \mathcal I \) only if a charging station is opened there, and that the number of installed chargers at each location $i$ is bounded above by $M_i$.
For each scenario $\omega\in\Omega$, the recourse function is defined as
\begin{align}
Q^\omega(z)=
\min_{x,u}
&\quad 
\sum_{i\in\mathcal I}\sum_{j\in\mathcal J} c_{ij}x_{ij}
+
\sum_{j\in\mathcal J}p_j u_{j}
-
\sum_{j\in\mathcal J}r_jd^{\omega}_{j}\\
\text{s.t.}& \quad \sum_{i\in\mathcal I}x_{ij}+u_{j}=d^{\omega}_{j},\ \forall j\in\mathcal J \label{eq:cons3}\\
&\quad \sum_{j\in\mathcal J}x_{ij}\le C_i z_i,\ \forall i\in\mathcal I \label{eq:cons4}\\
&\quad x_{ij}\ge 0, \ u_{j}\ge0, \ \forall i\in\mathcal I,\ j\in\mathcal J \label{eq:cons5}
\end{align}
where constraints \eqref{eq:cons3} ensure that the demand at each customer site $j\in \mathcal J$ is either satisfied by an open charging station or accounted for as unmet demand. Constraints \eqref{eq:cons4} enforce capacity limits by requiring that the total demand served at each open station $i\in \mathcal I$ does not exceed the product of the number of installed chargers $z_i$ and the per-charger capacity $C_i$. Finally, constraints \eqref{eq:cons5} impose non-negativity on both satisfied and unmet demand variables.

In our experiments, we set $|\mathcal I|=8$ and $|\mathcal J|=12$. The other parameters are uniformly sampled within the following ranges:
$f_i \in [80,300]$, $b_i \in [5,40]$, $c_{ij} \in [1,80]$, 
$p_j \in [30,120]$, $r_j \in [5,60]$, 
$C_i \in [40,120]$, and $M_i \in [10,80]$.

 

\section{Additional Numerical Results}\label{append:results}
\begin{figure*}[htbp]
    \centering

    \begin{subfigure}[b]{0.32\textwidth}
        \centering
        \includegraphics[width=\textwidth]{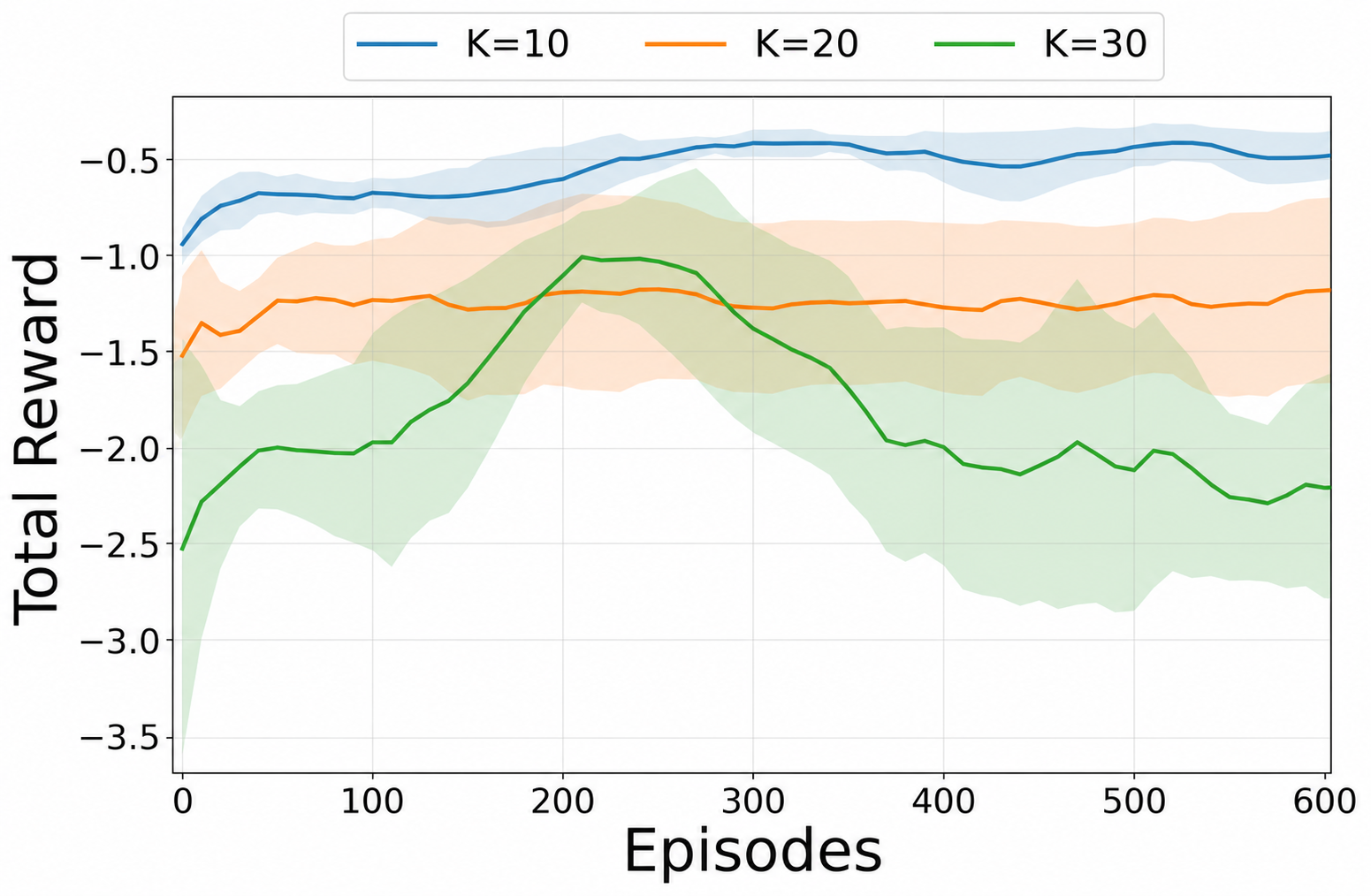}
        \caption{Total Reward}
        \label{fig:k_total_reward}
    \end{subfigure}
    \hfill
    \begin{subfigure}[b]{0.32\textwidth}
        \centering
        \includegraphics[width=\textwidth]{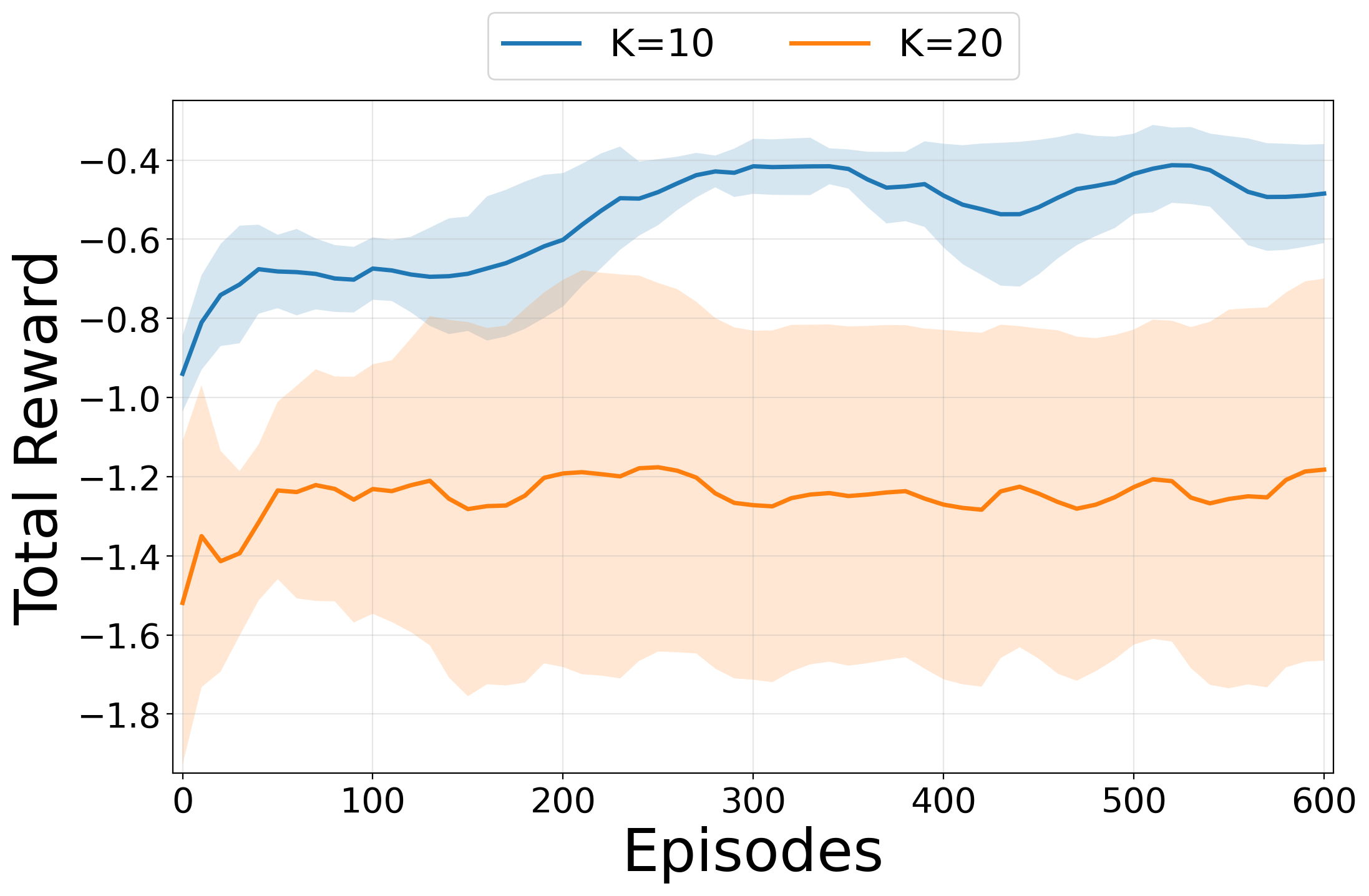}
        \caption{Total Reward (Close-up)}
        \label{fig:k_total_reward_zoom}
    \end{subfigure}
    \hfill
    \begin{subfigure}[b]{0.32\textwidth}
        \centering
        \includegraphics[width=\textwidth]{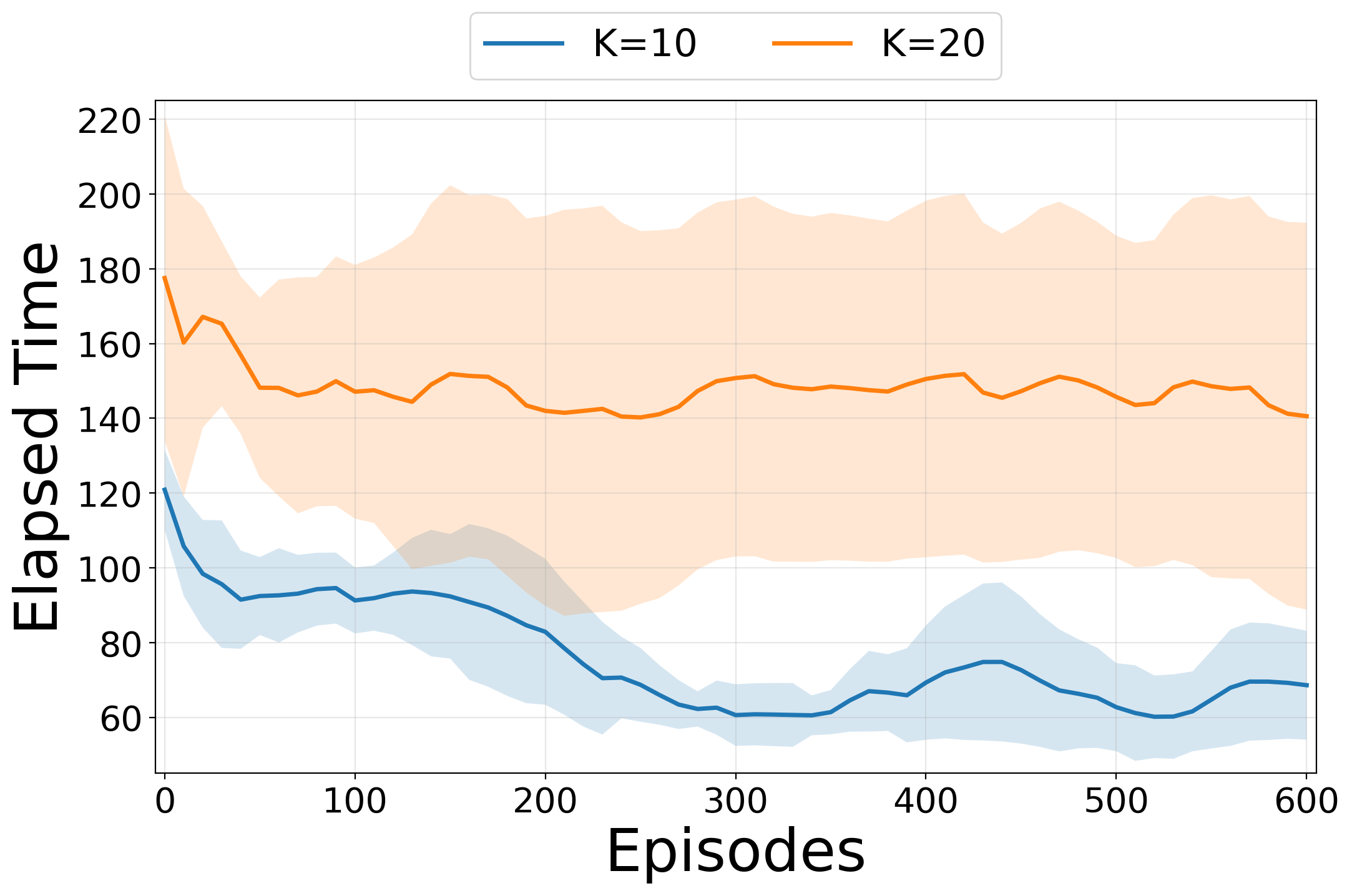}
        \caption{Elapsed Time (Close-up)}
        \label{fig:k_time_zoom}
    \end{subfigure}

    \caption{Evaluation results under different numbers of cuts added per iteration $K$.}
    \label{fig:k_tuning_results_append}
\end{figure*}

\end{document}